\documentclass[12pt,reqno]{amsart}
\usepackage[T1]{fontenc}
\usepackage{lmodern}
\usepackage{microtype}
\usepackage[margin=28mm]{geometry}
\usepackage{amsmath,amssymb,amsthm,mathtools}
\usepackage{xcolor}
\usepackage{enumitem}
\usepackage{hyperref}
\hypersetup{colorlinks=true,linkcolor=blue!45!black,citecolor=blue!45!black,urlcolor=blue!45!black}
\numberwithin{equation}{section}
\newtheorem{theorem}{Theorem}[section]
\newtheorem{conjecture}{Conjecture}
\newtheorem{lemma}[theorem]{Lemma}
\newtheorem{corollary}[theorem]{Corollary}
\newtheorem{proposition}[theorem]{Proposition}
\DeclareMathOperator{\Ric}{Ric}
\DeclareMathOperator{\Hess}{Hess}
\newcommand{\Sph}{\mathbb S}
\newcommand{\RP}{\mathbb{RP}}
\newcommand{\R}{\mathbb R}

\title[Scalar curvature deformation of isoparametric tubes]{Scalar curvature deformation of isoparametric tubes}

\author[H. X. Tan]{Huixin Tan}
\address{School of Mathematical Sciences, Laboratory of Mathematics and Complex Systems, Beijing Normal University, Beijing 100875, P. R.
CHINA.}
\email{hxtan@mail.bnu.edu.cn}

\author[Y. H. Zhang]{Yunheng Zhang$^{*}$}
\address{School of Mathematical Sciences, Laboratory of Mathematics and Complex Systems, Beijing Normal University, Beijing 100875, P. R.
CHINA.}
\email{yunheng@mail.bnu.edu.cn}

\thanks{$^{*}$ Corresponding author.}
\begin{document}

\begin{abstract}
	We resolve Ge's local scalar curvature rigidity conjecture for spherical
	isoparametric tubes.
	The conjecture fails throughout its stated radius range for families
	with at least two distinct principal curvatures, via arbitrarily small
	smooth deformations supported in any regular annulus.
	Geodesic balls of radius less than \(\frac{\pi}{2}\) are locally rigid
	under \(C^2\) smallness.
	We also prescribe sufficiently small smooth nonradial scalar-curvature
	increments with fixed compact interior support, under antipodal symmetry
	for families with 2, 4, or 6 distinct principal curvatures and
	without symmetry for families with 3 distinct principal curvatures.
\end{abstract}

\maketitle

\section{Introduction}

We study Ge's conjecture \cite[Problem~4, Conjecture~1]{Ge2019} on local scalar curvature rigidity for isoparametric tubes in the unit sphere \(\Sph^{n+1}\), where \(n\ge2\).

Scalar-curvature boundary rigidity is motivated in part by the positive mass theorem and its geometric consequences.
On the Euclidean unit ball, a metric with nonnegative scalar curvature, the standard induced boundary metric and boundary mean curvature at least its Euclidean value must be isometric to the Euclidean metric; see Miao \cite{Miao2002} and \cite[Theorem~1.4]{Brendle2012}.
Shi and Tam \cite{ShiTam2002} established mean-curvature comparison and rigidity results for strictly convex Euclidean boundaries.

The analogous spherical problem led to Min-Oo's conjecture \cite{MinOo1998}: a metric on the hemisphere in \(\Sph^{n+1}\), with scalar curvature at least \(n(n+1)\), the standard induced boundary metric and totally geodesic boundary, should be round.
Brendle, Marques and Neves \cite{BrendleMarquesNeves2011} disproved this conjecture in dimensions at least three.
Their perturbative construction also rules out the corresponding local rigidity statement with a boundary mean-curvature lower bound \cite[Theorem~4]{BrendleMarquesNeves2011}; see \cite[the remark following Theorem~3]{BrendleMarques2011}.
They further constructed examples agreeing with the round metric on a neighborhood of the boundary \cite[Theorem~7]{BrendleMarquesNeves2011}.

On smaller domains, Brendle and Marques \cite[Theorem~3]{BrendleMarques2011} proved local rigidity for geodesic balls \(B_\theta\subset\Sph^{n+1}\) with \(\cos\theta\ge\frac{2}{\sqrt{n+4}}\).
For sufficiently small \(C^2\) perturbations of the round metric, scalar curvature at least \(n(n+1)\), preservation of the induced boundary metric and boundary mean curvature at least its round value imply isometry by a diffeomorphism fixing the boundary pointwise.
Cox, Miao and Tam \cite{CoxMiaoTam2013} subsequently enlarged this radius range.
These results led Ge to ask whether local rigidity holds on every geodesic ball of radius less than \(\frac{\pi}{2}\) and whether it extends to tubes in other isoparametric families.

For a Cartan--M\"unzner polynomial \(F\) of degree \(g\in\{1,2,3,4,6\}\), normalized so that \(f=F|_{\Sph^{n+1}}\) has range \([-1,1]\), write \(M_t=f^{-1}(t)\) and \(M_\pm=f^{-1}(\pm1)\).
The regular leaves have \(g\) distinct principal curvatures.
If \(m_\pm+1\) are the codimensions of \(M_\pm\), then \(M_{c_0}\) is the unique minimal regular leaf, where \(c_0=\frac{m_--m_+}{m_-+m_+}\).
The tubes about \(M_+\) are
\begin{equation}\label{eq:Ge-tube}
	\Sph^{n+1}_{F+c}=\{f\ge c\}=\{r\le\theta\},
	\qquad \theta=\frac1g\arccos c,\quad -1<c<1,
\end{equation}
where \(r\) is the round distance to \(M_+\) and \(f=\cos(gr)\) on the regular region.
The regular leaves and focal submanifolds are connected by M\"unzner's structure theorems \cite{Munzner1980,Munzner1981}; see also \cite[Section~3.4]{CecilRyan2025}.
Thus each tube in \eqref{eq:Ge-tube} is a compact connected domain with smooth boundary \(M_c\).

For a Riemannian metric \(\mathbf g\), let \(R_{\mathbf g}\) denote its scalar curvature and \(\mathbf g^T\) the induced boundary metric.
Let \(\nu_{\mathrm{in}}\) and \(\nu_{\mathrm{out}}\) denote the inward and outward unit normals to the boundary with respect to \(\mathbf g\), respectively.
We use the second fundamental form \(\mathrm{II}_{\mathbf g}(X,Y)=-\mathbf g(\nabla_X\nu_{\mathrm{in}},Y)=\mathbf g(\nabla_X\nu_{\mathrm{out}},Y)\) and mean curvature \(H_{\mathbf g}=\operatorname{tr}\mathrm{II}_{\mathbf g}\).
The Laplacian is \(\Delta_{\mathbf g}=\operatorname{tr}_{\mathbf g}\Hess_{\mathbf g}\), and the curvature convention is \(\mathcal R(X,Y)Z=\nabla_X\nabla_YZ-\nabla_Y\nabla_XZ-\nabla_{[X,Y]}Z\).
Let \(\bar{\mathbf g}\) denote the round metric on \(\Sph^{n+1}\).
Unless specified otherwise, norms are computed with respect to \(\bar{\mathbf g}\).
We write \(d\mu\) and \(d\sigma\) for the volume and induced boundary area measures of \(\bar{\mathbf g}\), respectively.

Ge \cite{Ge2019} proposed extending local scalar curvature rigidity from geodesic balls to isoparametric tubes.
When \(g=1\), these tubes are geodesic balls and \(c_0=0\).
The proposed range \(c_0<c<1\), equivalently \(0<\theta<\frac1g\arccos c_0\), stops before the boundary reaches the minimal regular leaf; in degree \(1\), its excluded endpoint is precisely the hemisphere.
Thus the conjecture asks both for rigidity throughout the open geodesic-ball range and for its extension to the other isoparametric families.
Ge retained the \(C^2\)-smallness hypothesis and also considered replacing the mean-curvature bound by the stronger comparison of second fundamental forms.

\begin{conjecture}[Ge \cite{Ge2019}]\label{conj:Ge}
	Let \(\mathbf g\) be a Riemannian metric on the isoparametric tube \(\Sph^{n+1}_{F+c}\) defined in \eqref{eq:Ge-tube}, for \(c\in(c_0,1)\), with the following properties:
	\begin{itemize}[leftmargin=2em]
		\item \(R_{\mathbf g}\ge n(n+1)\) at each point in \(\Sph^{n+1}_{F+c}\).
		\item The metrics \(\mathbf g\) and \(\bar{\mathbf g}\) induce the same metric on \(M_c=\partial\Sph^{n+1}_{F+c}\).
		\item \(H_{\mathbf g}\ge H_{\bar{\mathbf g}}\) at each point on \(M_c\).
	\end{itemize}
	If \(\mathbf g-\bar{\mathbf g}\) is sufficiently small in the \(C^2\)-norm, then \(\varphi^*\mathbf g=\bar{\mathbf g}\) for some diffeomorphism \(\varphi\colon\Sph^{n+1}_{F+c}\to\Sph^{n+1}_{F+c}\) with \(\varphi|_{M_c}=\mathrm{id}\).
\end{conjecture}

Ge also proposed a global rigidity conjecture under the Ricci-curvature lower bound \(\Ric_{\mathbf g}\ge n\mathbf g\), without requiring \(\mathbf g\) to be close to the round metric \(\bar{\mathbf g}\) \cite[Problem~4, Conjecture~2]{Ge2019}.
For \(g=2\), Lai \cite[Corollary~1.1]{Lai2020} proved this global rigidity for the corresponding Clifford tubes under weaker boundary curvature assumptions.

We prove that Conjecture~\ref{conj:Ge} holds for geodesic balls and fails for every family of higher degree.
For the latter families, we prescribe small radial scalar-curvature increments by deformations supported in an arbitrary regular annulus.

\begin{theorem}\label{thm:radial}
	Suppose \(g\in\{2,3,4,6\}\), \(-1<c<1\), and \(0<\alpha<\beta<\theta=\frac{1}{g}\arccos c\).
	For every \(\chi\in C_c^\infty((\alpha,\beta))\), extended by zero to \([0,\theta]\), there exist \(\varepsilon_0>0\) and a smooth family of metrics \(\mathbf{g}_{\varepsilon}\), \(|\varepsilon|<\varepsilon_0\), on the fixed tube \(\Sph^{n+1}_{F+c}\) with the following properties:
	\begin{enumerate}[label=\textup{(\roman*)},leftmargin=2.2em]
		\item \(\mathbf{g}_0=\bar{\mathbf{g}}\), and \(\mathbf{g}_{\varepsilon}\to\bar{\mathbf{g}}\) in \(C^k\) for every nonnegative integer \(k\) as \(\varepsilon\to0\).
		\item \(\mathbf{g}_{\varepsilon}=\bar{\mathbf{g}}\) outside \(\{\alpha<r<\beta\}\).
		In particular, the metrics agree with the round metric on fixed neighborhoods of \(M_+\) and \(M_c=\partial\Sph^{n+1}_{F+c}\).
		\item \(R_{\mathbf{g}_{\varepsilon}}=n(n+1)+\varepsilon\chi(r)\).
		\item If \(\chi\ge0\) is nonzero and \(0<\varepsilon<\varepsilon_0\), there is no diffeomorphism \(\varphi\colon\Sph^{n+1}_{F+c}\to\Sph^{n+1}_{F+c}\) such that \(\varphi^*\mathbf{g}_{\varepsilon}=\bar{\mathbf{g}}\).
	\end{enumerate}
	Consequently, for every isoparametric family of degree \(g\in\{2,3,4,6\}\), Conjecture~\ref{conj:Ge} fails on every tube \(\Sph^{n+1}_{F+c}\) in its stated range \(c_0<c<1\).
\end{theorem}

The metrics in Theorem~\ref{thm:radial} agree with \(\bar{\mathbf g}\) near \(M_c\), so \(\mathrm{II}_{\mathbf g_\varepsilon}=\mathrm{II}_{\bar{\mathbf g}}\) on \(M_c\).
Thus, for \(g\ge2\), Conjecture~\ref{conj:Ge} remains false under the stronger boundary condition \(\mathrm{II}_{\mathbf g}\ge\mathrm{II}_{\bar{\mathbf g}}\) considered by Ge.
In comparison with the hemisphere examples of \cite[Theorem~7]{BrendleMarquesNeves2011}, Theorem~\ref{thm:radial} gives nonrigidity with unchanged boundary geometry at every radius in Ge's conjectured range, including arbitrarily small tubes about \(M_+\), and with arbitrarily small smooth perturbations.
Akutagawa's non-Einstein relative Yamabe metric on the solid torus bounded by the minimal Clifford torus \cite[Section~2]{Akutagawa2021} satisfies a different boundary condition at \(c=c_0\).

For geodesic balls, Gong and Wei \cite[Theorem~1.1]{GongWei2026} proved local rigidity throughout \(0<\theta<\frac{\pi}{2}\) under \(C^{2,\gamma}\) smallness, where \(0<\gamma<1\).
Our treatment of this case uses their positivity theorem together with a \(W^{2,p}\) remainder estimate, for \(p>n+1\), to obtain rigidity under the \(C^2\)-smallness hypothesis in Ge's original formulation; see Section~\ref{sec:balls}.

\begin{theorem}\label{thm:balls}
	Let \(n\ge2\) and \(0<\theta<\frac{\pi}{2}\).
	Every smooth metric \(\mathbf{g}\) on the round geodesic ball \(B_\theta\subset\Sph^{n+1}\) sufficiently close to \(\bar{\mathbf{g}}\) in \(C^2\), with \(R_{\mathbf{g}}\ge n(n+1)\), \(\mathbf{g}^T=\bar{\mathbf{g}}^T\), and \(H_{\mathbf{g}}\ge H_{\bar{\mathbf{g}}}\), satisfies \(\varphi^*\mathbf{g}=\bar{\mathbf{g}}\) for a smooth diffeomorphism \(\varphi\) fixing the boundary pointwise.
\end{theorem}

Theorem~\ref{thm:balls} extends the \(C^2\) local rigidity theorem of Brendle and Marques \cite[Theorem~3]{BrendleMarques2011} to all geodesic balls of radius less than \(\frac{\pi}{2}\). Hence, Theorems~\ref{thm:radial} and~\ref{thm:balls} settle Conjecture~\ref{conj:Ge} throughout its stated range.

We also prescribe sufficiently small smooth nonradial scalar-curvature increments supported in a fixed compact subset of the interior: Corollary~\ref{cor:antipodal} treats tubes with \(g\in\{2,4,6\}\) under antipodal symmetry, while Proposition~\ref{prop:cubic-prescription} treats tubes with \(g=3\) without symmetry.
When \(g=3\), the metric can be chosen in any prescribed \(C^{2,\gamma}\) neighborhood of \(\bar{\mathbf g}\), with the \(C^{0,\gamma}\) smallness bound depending on that neighborhood; differentiable dependence at \(\bar{\mathbf g}\) is not asserted.
Both constructions use Corvino's local deformation theorem \cite{Corvino2000,CorvinoPollack2011}, applied on the antipodal quotient when \(g\in\{2,4,6\}\) and to nearby metrics of constant scalar curvature obtained from the radial construction when \(g=3\).

Section~\ref{sec:prescription} proves Theorem~\ref{thm:radial} by a radial scalar-curvature deformation.
Section~\ref{sec:balls} proves Theorem~\ref{thm:balls}.
Section~\ref{sec:nonradial} treats nonradial scalar-curvature prescription on tubes with \(g\in\{2,3,4,6\}\).

\section{Radial scalar-curvature deformations}\label{sec:prescription}

We prove Theorem~\ref{thm:radial} by prescribing the scalar curvature on a regular annulus.
The same construction provides the metrics of constant scalar curvature used in Section~\ref{sec:nonradial}.

\begin{proof}[Proof of Theorem~\ref{thm:radial}]
	The normal flow identifies \(\{\alpha\le r\le\beta\}\) with \([\alpha,\beta]\times M\), where \(M\) is a regular leaf.
	Write \(\bar{\mathbf g}=dr^2+\mathbf h_r\), and let \(S=\frac12\mathbf h_r^{-1}\mathbf h_r'\) and \(H=\operatorname{tr}S\) be the shape operator and mean curvature of the regular leaves with respect to the unit normal \(\partial_r\).
	Primes denote differentiation with respect to \(r\).

	For tangent vectors \(X,Y\) to a regular leaf, the Gauss equation and the normal Riccati equation take the form
	\begin{equation}\label{eq:gauss-riccati}
		\begin{aligned}
			\mathbf h_r(\mathcal R_{\mathbf h_r}(X,Y)Y,X)
			&=\bar{\mathbf g}(\mathcal R_{\bar{\mathbf g}}(X,Y)Y,X)
			+\mathbf h_r(SX,X)\mathbf h_r(SY,Y)
			-\mathbf h_r(SX,Y)^2,\\
			S'+S^2+\mathcal R_{\bar{\mathbf g}}(\,\cdot\,,\partial_r)\partial_r&=0.
		\end{aligned}
	\end{equation}
	Since \(\bar{\mathbf g}\) has constant sectional curvature one, taking traces in~\eqref{eq:gauss-riccati} gives
	\begin{equation}\label{eq:traced-gauss-riccati}
		R_{\mathbf h_r}=n(n-1)+H^2-|S|^2, \qquad
		H'=-|S|^2-n.
	\end{equation}

	For \(N=N(r)>0\), put \(\widehat{\mathbf g}_N=N^2dr^2+\mathbf h_r\) and \(R_N=R_{\widehat{\mathbf g}_N}\).
	The unit normal to the regular leaves with respect to \(\widehat{\mathbf g}_N\) is \(N^{-1}\partial_r\), and the corresponding shape operator and mean curvature are \(\frac{1}{N}S\) and \(\frac{1}{N}H\), respectively.
	Since \(N\) depends only on \(r\), the integral curves of \(N^{-1}\partial_r\) are geodesics.

	Applying the traced Gauss equation to the leaves in the metric \(\widehat{\mathbf g}_N\) gives
	\begin{equation}\label{Guasstr1}
			R_{\mathbf h_r}
		=R_N-2\Ric_{\widehat{\mathbf g}_N}
		(N^{-1}\partial_r,N^{-1}\partial_r)
		+\frac{H^2-|S|^2}{N^2}.
	\end{equation}
	Tracing the normal Riccati equation along the unit normal \(N^{-1}\partial_r\) gives
	\[
		\Ric_{\widehat{\mathbf g}_N}
		(N^{-1}\partial_r,N^{-1}\partial_r)
		=-\frac1N\left(\frac HN\right)'
		-\frac{|S|^2}{N^2}.
	\]
	Substituting this identity into \eqref{Guasstr1}, we obtain
	\begin{equation}\label{eq:lapse-scalar}
		\begin{aligned}
			R_N
			&=R_{\mathbf h_r}
			-\frac{H^2+|S|^2}{N^2}
			-\frac2N\left(\frac HN\right)'\\
			&=R_{\mathbf h_r}
			-\frac{H^2+|S|^2+2H'}{N^2}
			+\frac{2HN'}{N^3}.
		\end{aligned}
	\end{equation}
	Since the principal curvatures are constant on each regular leaf, \(H\), \(|S|^2\), and hence \(R_N\) depend only on \(r\).
	For a radial function \(u=u(r)\), its gradient with respect to \(\widehat{\mathbf g}_N\) is \(\operatorname{grad}_{\widehat{\mathbf g}_N}u=N^{-2}u'\partial_r\).
	Then the divergence formula gives
	\begin{equation}\label{radial Laplacian}
			\begin{aligned}
			\Delta_{\widehat{\mathbf g}_N}u
			=\operatorname{div}_{\widehat{\mathbf g}_N}
			\left(\frac{u'}{N}N^{-1}\partial_r\right)
			=\frac1N\left(\frac{u'}N\right)'
			+\frac{Hu'}{N^2}
			=\frac1{N^2}
			\left[u''+\left(H-\frac{N'}N\right)u'\right].
		\end{aligned}
	\end{equation}
	For any smooth \(u>0\), set \(\omega=\frac{2}{n-1}\log u\), so that \(e^{2\omega}=u^{\frac{4}{n-1}}\).
	The conformal change formula in dimension \(n+1\) gives
	\[
		R_{e^{2\omega}\widehat{\mathbf g}_N}
		=e^{-2\omega}\bigl(R_N-2n\Delta\omega-n(n-1)|\nabla\omega|^2\bigr),
	\]
	where the derivatives and norms are taken with respect to \(\widehat{\mathbf g}_N\).
	By the chain rule,
	\[
		\Delta\omega=\frac{2}{n-1}\left(\frac{\Delta u}{u}-\frac{|\nabla u|^2}{u^2}\right),
		\qquad
		|\nabla\omega|^2=\frac{4}{(n-1)^2}\frac{|\nabla u|^2}{u^2}.
	\]
	The terms involving \(|\nabla u|^2\) cancel, leaving \(R_N-\frac{4n}{n-1}\frac{\Delta u}{u}\) in parentheses.
	Multiplying by \(e^{-2\omega}=u^{-\frac{4}{n-1}}\) and collecting powers of \(u\), we obtain
	\begin{equation}\label{conformal}
		R_{u^{\frac{4}{n-1}}\widehat{\mathbf g}_N}
		=u^{-\frac{n+3}{n-1}}
		\left(-\frac{4n}{n-1}\Delta_{\widehat{\mathbf g}_N}u+R_Nu\right).
	\end{equation}
	Now suppose \(u=u(r)\).
	Substituting the expression for \(\Delta_{\widehat{\mathbf g}_N}u\) from \eqref{radial Laplacian} into the conformal formula \eqref{conformal}, we find that \(u^{\frac{4}{n-1}}\widehat{\mathbf g}_N\) has scalar curvature \(n(n+1)+\varepsilon\chi\) if and only if
	\begin{equation}\label{eq:radial-ode}
		u''+\left(H-\frac{N'}N\right)u'
		=\frac{(n-1)N^2}{4n}\left[R_Nu-(n(n+1)+\varepsilon\chi)u^{\frac{n+3}{n-1}}\right].
	\end{equation}

	Fix \(q_0\in C_c^\infty((\alpha,\beta))\).
	For \(s\in\R\) and \(\boldsymbol\lambda=(\lambda_1,\lambda_2)\in\R^2\), set
	\begin{equation}\label{eq:lapse-parameters}
		N=1+sq_0+\lambda_1q_1+\lambda_2q_2,
	\end{equation}
	where \(q_1,q_2\in C_c^\infty((\alpha,\beta))\) will be chosen below.
	We solve \eqref{eq:radial-ode} with \(u(\alpha)=1\) and \(u'(\alpha)=0\).
	For each sufficiently small pair \((s,\varepsilon)\), we seek \(\lambda_1,\lambda_2\) such that \(u(\beta)=1\) and \(u'(\beta)=0\).
	At the parameter origin \((s,\boldsymbol\lambda,\varepsilon)=(0,0,0)\), \(N=u=1\).
	By smooth dependence for ordinary differential equations, the solution \(u\) exists on \([\alpha,\beta]\), remains positive, and depends smoothly on \((s,\boldsymbol\lambda,\varepsilon)\) for all sufficiently small parameters.

	To impose the two conditions at \(r=\beta\), define the endpoint map 
	\[\mathcal E(s,\boldsymbol\lambda,\varepsilon)=(u(\beta)-1,u'(\beta)).
	\]
	We seek \(q_1,q_2\) such that the derivative \(\partial_{\boldsymbol\lambda}\mathcal E(0,0,0)\) is invertible.
	To compute this derivative, set \(z_j=\left.\partial_{\lambda_j}u\right|_0\) for \(j=1,2\).
	By~\eqref{eq:lapse-parameters}, \(\partial_{\lambda_j}N=q_j\).
	Differentiating \eqref{eq:lapse-scalar} at the origin and using \eqref{eq:traced-gauss-riccati}, we obtain
	\[
		\begin{aligned}
			\left.\partial_{\lambda_j}R_N\right|_0
			&=2(H^2+|S|^2+2H')q_j+2Hq_j'\\
			&=2Hq_j'+2(H^2-|S|^2-2n)q_j.
		\end{aligned}
	\]
	At the origin, \(u'=0\) and the bracket on the right-hand side of \eqref{eq:radial-ode} vanishes, so differentiating \(\frac{N'}{N}\) and the prefactor \(N^2\) produces no additional first-order terms.
	Differentiating the remaining terms gives
	\[
		\begin{aligned}
			z_j''+Hz_j'
			&=\frac{n-1}{4n}\left[\left.\partial_{\lambda_j}R_N\right|_0
			+n(n+1)\left(1-\frac{n+3}{n-1}\right)z_j\right]\\
			&=\frac{n-1}{4n}\left.\partial_{\lambda_j}R_N\right|_0-(n+1)z_j.
		\end{aligned}
	\]
	Let \(P\) denote the operator obtained by moving the zeroth-order term to the left.
	Substituting the expression above for \(\left.\partial_{\lambda_j}R_N\right|_0\) and differentiating the initial conditions at \(r=\alpha\), we obtain
	\begin{equation}\label{eq:radial-linearization}
		\begin{aligned}
			Pz&=z''+Hz'+(n+1)z,\\
			\frac{4n}{n-1}Pz_j&=2Hq_j'+2(H^2-|S|^2-2n)q_j,\qquad z_j(\alpha)=z_j'(\alpha)=0.
		\end{aligned}
	\end{equation}
	Let \(y_1,y_2\) solve \(Py_i=0\), with \((y_1,y_1')(\beta)=(1,0)\) and \((y_2,y_2')(\beta)=(0,1)\), and set \(k_i=Hy_i'+ny_i\).
	The functions \(k_i\) arise when the \(q_j'\) term in \eqref{eq:radial-linearization} is integrated by parts.
	We first show that \(k_1\) and \(k_2\) are linearly independent on every nonempty open subinterval of \((\alpha,\beta)\).

	Fix \(b_1,b_2\in\R\), and set \(y=b_1y_1+b_2y_2\) and \(k=b_1k_1+b_2k_2\).
	Then \(Py=0\) and \(k=Hy'+ny\).
	Using \eqref{eq:traced-gauss-riccati} and \(Py=0\), we obtain
	\[
		k'=(H'+n)y'+Hy''=-(|S|^2+H^2)y'-(n+1)Hy.
	\]
	Substituting \(ny=k-Hy'\) gives
	\begin{equation}\label{eq:coupling}
		k'=-\left(|S|^2-\frac{H^2}{n}\right)y'-\frac{(n+1)H}{n}k.
	\end{equation}
	Since \(g\ge2\), the regular leaves have at least two distinct principal curvatures, and hence \(|S|^2-\frac{H^2}{n}=|S-\frac{H}{n}I|^2>0\).

	Suppose that \(k=0\) on a nonempty open subinterval of \((\alpha,\beta)\).
	Then \(k'=0\) there, so \eqref{eq:coupling} and the strict inequality above imply \(y'=0\).
	The identity \(k=Hy'+ny\) then gives \(y=0\) on the same subinterval.
	By uniqueness for \(Py=0\), we have \(y\equiv0\) on \([\alpha,\beta]\).
	Evaluating \(y\) and \(y'\) at \(r=\beta\), we obtain \(b_1=y(\beta)=0\) and \(b_2=y'(\beta)=0\).
	Thus \(k_1,k_2\) are linearly independent on every nonempty open subinterval of \((\alpha,\beta)\).

	Choose a nonzero nonnegative \(\eta\in C_c^\infty((\alpha,\beta))\) and set \(q_j=\eta k_j\).
	Let \(w>0\) satisfy \(w'=Hw\), so that \(P=w^{-1}\partial_r(w\partial_r)+(n+1)\).
	Using \eqref{eq:traced-gauss-riccati} and \(w'=Hw\), we obtain \((wHy_i)'-w(H^2-|S|^2-2n)y_i=wk_i\).
	Multiplying the equation for \(z_j\) in~\eqref{eq:radial-linearization} by \(wy_i\) and integrating by parts therefore gives
	\[
		\frac{4n}{n-1}\int_\alpha^\beta wy_iPz_j\,dr=-2\int_\alpha^\beta wk_iq_j\,dr.
	\]
	There are no boundary terms because \(q_j\) is compactly supported.
	Green's identity yields
	\begin{equation}\label{eq:endpoint-identity}
		\begin{aligned}
			w(\beta)\bigl[y_i(\beta)z_j'(\beta)-y_i'(\beta)z_j(\beta)\bigr]
			&=\bigl[w(y_i z_j'-y_i'z_j)\bigr]_\alpha^\beta
			=\int_\alpha^\beta wy_iPz_j\,dr\\
			&=-\frac{n-1}{2n}\int_\alpha^\beta wk_iq_j\,dr
			=-\frac{n-1}{2n}\int_\alpha^\beta w\eta k_i k_j\,dr.
		\end{aligned}
	\end{equation}
	The \(j\)th column of \(\partial_{\boldsymbol\lambda}\mathcal E(0,0,0)\) is \((z_j(\beta),z_j'(\beta))^{\mathsf T}\).
	Substituting the terminal data for \(y_1,y_2\) into \eqref{eq:endpoint-identity} therefore gives
	\[
		w(\beta)\begin{pmatrix}0&1\\-1&0\end{pmatrix}
		\partial_{\boldsymbol\lambda}\mathcal E(0,0,0)
		=-\frac{n-1}{2n}\left(\int_\alpha^\beta w\eta\,k_i k_j\,dr\right)_{i,j=1}^2.
	\]
	The matrix of integrals is positive definite because \(w\eta\ge0\) and \(k_1,k_2\) are linearly independent on an interval where \(w\eta>0\).
	Since \(n\ge2\), \(w(\beta)>0\), and the displayed constant matrix has determinant one, \(\partial_{\boldsymbol\lambda}\mathcal E(0,0,0)\) is invertible.

	The implicit function theorem gives smooth parameters \(\boldsymbol\lambda(s,\varepsilon)\), vanishing at the origin, such that \(\mathcal E(s,\boldsymbol\lambda(s,\varepsilon),\varepsilon)=0\).
	For these parameters, define
	\begin{equation}\label{eq:radial-family}
		\mathbf g_{s,\varepsilon}=u^{\frac4{n-1}}(N^2dr^2+\mathbf h_r),
		\qquad R_{\mathbf g_{s,\varepsilon}}=n(n+1)+\varepsilon\chi(r),
	\end{equation}
	and extend the metric by \(\bar{\mathbf g}\) outside the annulus.
	All the functions \(q_0,q_1,q_2,\chi\) vanish on fixed neighborhoods of both endpoints.
	There \(N=1\), and uniqueness for \eqref{eq:radial-ode}, applied separately at the two endpoints where \(u=1\) and \(u'=0\), yields \(u=1\) on these neighborhoods.
	The extension is smooth, including across the focal submanifold, and depends smoothly on \((s,\varepsilon)\).
	Compactness gives \(\mathbf g_{s,\varepsilon}\to\bar{\mathbf g}\) in every \(C^k\) norm as \((s,\varepsilon)\to(0,0)\).
	The required family is \(\mathbf g_\varepsilon=\mathbf g_{0,\varepsilon}\).

	Finally, if \(\chi\ge0\) is nonzero and \(\varepsilon>0\), then \(R_{\mathbf g_\varepsilon}\ge n(n+1)\), with strict inequality where \(\chi>0\).
	An isometry to \(\bar{\mathbf g}\) would force \(R_{\mathbf g_\varepsilon}\equiv n(n+1)\), since scalar curvature is invariant under pullback.
	The metric is round near \(M_c\), so its induced boundary metric and second fundamental form agree with those of \(\bar{\mathbf g}\).
	For \(c_0<c<1\), this gives the asserted counterexamples to Conjecture~\ref{conj:Ge}.
\end{proof}

We also record the case of constant scalar curvature for use in Proposition~\ref{prop:cubic-prescription}.
Taking \(\varepsilon=0\) in \eqref{eq:radial-family} gives nearby metrics \(\mathbf g_{s,0}\) with scalar curvature \(n(n+1)\).
The functions \(q_1,q_2\) above do not depend on \(q_0\).
Since their supports lie in \(\operatorname{supp}\eta\Subset(\alpha,\beta)\), we may choose a nonzero \(q_0\in C_c^\infty((\alpha,\beta))\) with support disjoint from them and apply the same construction with this choice.
On \(\operatorname{supp}q_0\), \eqref{eq:lapse-parameters} then reduces to \(N-1=sq_0\), so \(N\not\equiv1\) for \(s\ne0\).
Consequently, \(\mathbf g_{s,0}\ne\bar{\mathbf g}\): equality would force \(u=1\) from the tangential components and then \(N=1\) from the radial component.

\section{Rigidity of geodesic balls}\label{sec:balls}

We now prove the rigidity statement in Theorem~\ref{thm:balls}.
The key ingredient is Gong and Wei's positivity theorem~\cite[Theorem~4.3]{GongWei2026}.
We follow the compactness argument in~\cite[Section~5.2]{GongWei2026} and supply a \(W^{2,p}\) estimate for the Taylor remainder, which yields the \(C^2\) smallness statement required here.
In fact, the proof gives rigidity under \(W^{2,p}\) smallness for every fixed \(p>n+1\).
All norms, derivatives, and measures in this section are taken with respect to \(\bar{\mathbf{g}}\).

For a symmetric two-tensor \(h\), define the linearization of the scalar curvature at \(\mathbf g\) by \(L_{\mathbf g}h=\left.\frac{d}{dt}\right|_{t=0}R_{\mathbf g+th}\).
Explicitly,
\begin{equation}\label{def:L}
	L_{\mathbf g}h
	=-\Delta_{\mathbf g}(\operatorname{tr}_{\mathbf g}h)
	+\operatorname{div}_{\mathbf g}\operatorname{div}_{\mathbf g}h
	-\langle\Ric_{\mathbf g},h\rangle_{\mathbf g}.
\end{equation}
Its formal \(L^2\)-adjoint is
\begin{equation}\label{eq:adjoint}
	L_{\mathbf g}^*v=\Hess_{\mathbf g}v-(\Delta_{\mathbf g}v)\mathbf g-v\Ric_{\mathbf g}.
\end{equation}

To use the Brendle--Marques quadratic form~\cite{BrendleMarques2011} in the rigidity argument, we first make the metric perturbation divergence-free with respect to \(\bar{\mathbf g}\).
The following lemma is a consequence of Proposition~11 in \cite{BrendleMarques2011} and elliptic boundary regularity.

\begin{lemma}\label{lem:ball-gauge}
	Let \(0<\theta<\frac{\pi}{2}\) and \(p>n+1\).
	For every smooth Riemannian metric \(\mathbf g\) on the round geodesic ball \(B_\theta\) sufficiently close to \(\bar{\mathbf g}\) in \(W^{2,p}\), there is a smooth diffeomorphism \(\varphi\colon B_\theta\to B_\theta\), with \(\varphi|_{\partial B_\theta}=\mathrm{id}\), such that the tensor
	\begin{equation}\label{eq:ball-perturbation}
		h=\varphi^*\mathbf g-\bar{\mathbf g}
	\end{equation}
	satisfies
	\[
		\operatorname{div}_{\bar{\mathbf g}}h=0,\qquad
		\|h\|_{W^{2,p}}\le C\|\mathbf g-\bar{\mathbf g}\|_{W^{2,p}},
	\]
	where \(C\) depends only on \(n,\theta,p\).
	If \(\mathbf g^T=\bar{\mathbf g}^T\), then \(h^T=0\).
\end{lemma}

\begin{proof}
	The proof of Proposition~11 in \cite{BrendleMarques2011} yields a \(W^{3,p}\) diffeomorphism close to the identity and the stated estimate.
	The equation for \(\varphi\), with Dirichlet boundary data, is
	\[
		\operatorname{div}_{\bar{\mathbf g}}(\varphi^*\mathbf g-\bar{\mathbf g})=0,
		\qquad \varphi|_{\partial B_\theta}=\mathrm{id}.
	\]
	At \((\mathbf g,\varphi)=(\bar{\mathbf g},\mathrm{id})\), its linearization in \(\varphi\) is
	\(X\mapsto\operatorname{div}_{\bar{\mathbf g}}(\mathcal L_X\bar{\mathbf g})\), with \(X|_{\partial B_\theta}=0\).
	Up to an overall sign, the principal symbol is
	\(X\mapsto|\xi|^2X+\xi(X)\xi^\sharp\), which is positive definite for every nonzero covector \(\xi\).
	Here \(\xi^\sharp\) and the norm are defined by \(\bar{\mathbf g}\).
	Since \(p>n+1\), Sobolev embedding and the smallness of \(\mathbf g-\bar{\mathbf g}\) and \(\varphi-\mathrm{id}\) ensure that strong ellipticity persists at this solution.
	For smooth \(\mathbf g\), elliptic boundary regularity therefore bootstraps \(\varphi\) to smoothness up to \(\partial B_\theta\).
	Since \(\varphi\) fixes the boundary pointwise, \(\mathbf g^T=\bar{\mathbf g}^T\) implies \(h^T=0\).
\end{proof}

We now combine Lemma~\ref{lem:ball-gauge} with positivity of the Brendle--Marques quadratic form on the linearized constraint space and a \(W^{2,p}\) remainder estimate to prove Theorem~\ref{thm:balls}.

\begin{proof}[Proof of Theorem~\ref{thm:balls}]
	Let \(r\) be the round distance from the center, and set \(V=\cos r\) and \(c=\cos\theta>0\).
	As a height function on the unit sphere, \(V\) satisfies \(\Hess_{\bar{\mathbf g}}V=-V\bar{\mathbf g}\), and taking the trace gives \(\Delta_{\bar{\mathbf g}}V=-(n+1)V\).
	Since \(\Ric_{\bar{\mathbf g}}=n\bar{\mathbf g}\), substitution into~\eqref{eq:adjoint} yields
	\[
		L_{\bar{\mathbf g}}^*V=-V\bar{\mathbf g}+(n+1)V\bar{\mathbf g}-nV\bar{\mathbf g}=0.
	\]
	Moreover, \(V|_{\partial B_\theta}=\cos\theta=c\).
	Fix \(p>n+1\), and choose \(\varphi\) as in Lemma~\ref{lem:ball-gauge}, with \(h\) defined by~\eqref{eq:ball-perturbation}.
	Lemma~\ref{lem:ball-gauge} gives \(\operatorname{div}_{\bar{\mathbf g}}h=0\), \(h^T=0\) and the required \(W^{2,p}\) control.
	The curvature and boundary conditions are preserved by this pullback, so it suffices to prove that \(h=0\) when \(\|h\|_{W^{2,p}}\) is sufficiently small.

	Let \(\nu=\partial_r\) be the outward unit normal to \(\partial B_\theta\).
	For this perturbation \(h\), consider the Brendle--Marques quadratic form
	\begin{equation}\label{eq:ball-form}
		\begin{aligned}
			Q(h)={}&\frac{1}{4}\int_{B_\theta} V\bigl(|\nabla h|^2+|\nabla\operatorname{tr}h|^2+2|h|^2+2(\operatorname{tr}h)^2\bigr)\,d\mu\\
			&+\frac{(n+4)c^2-4}{4\sin\theta}\int_{\partial B_\theta} h(\nu,\nu)^2\,d\sigma
			+\frac{(n+2)c^2-1}{2\sin\theta}\int_{\partial B_\theta}|h(\nu,\cdot)^T|^2\,d\sigma.
		\end{aligned}
	\end{equation}
	The two boundary coefficients follow from \(V=c\), \(\partial_\nu V=-\sin\theta\), and \(H_{\bar{\mathbf g}}=n\cot\theta\) on \(\partial B_\theta\):
	\[
		\begin{aligned}
			\partial_\nu V+\frac{1}{4}H_{\bar{\mathbf g}}V
			&=\frac{(n+4)c^2-4}{4\sin\theta},\\
			\frac{1}{2}\partial_\nu V+\frac{n+1}{2n}H_{\bar{\mathbf g}}V
			&=\frac{(n+2)c^2-1}{2\sin\theta}.
		\end{aligned}
	\]
	Formula~\eqref{eq:ball-form} also defines a continuous quadratic form on \(H^1(B_\theta;S^2T^*B_\theta)\), the Sobolev space of symmetric covariant two-tensors.
	Define
	\begin{equation}\label{eq:ball-constraint-space}
		\begin{aligned}
			&Z=\bigl\{h\in H^1(B_\theta;S^2T^*B_\theta):
			\operatorname{div}_{\bar{\mathbf g}}h=0,\quad h^T=0,\\
			&\qquad\int_{B_\theta}\langle L_{\bar{\mathbf g}}^*\zeta,h\rangle\,d\mu=0
			\quad\text{for all }\zeta\in C^\infty(\overline{B_\theta})\bigr\}.
		\end{aligned}
	\end{equation}
	The divergence condition is understood in the distributional sense and the boundary condition in the trace sense.
	For smooth \(h\) with \(\operatorname{div}_{\bar{\mathbf g}}h=0\) and \(h^T=0\), Green's identity identifies the integral condition in the definition of \(Z\) with \(L_{\bar{\mathbf g}}h=0\) in \(B_\theta\) and vanishing first variation of the boundary mean curvature.
	This integral formulation imposes both linearized curvature constraints on \(H^1\) tensors without requiring a trace of their normal derivatives.
	The space \(Z\) is a closed linear subspace of \(H^1(B_\theta;S^2T^*B_\theta)\).
	Gong and Wei~\cite[Theorem~4.3]{GongWei2026} prove \(Q(h)>0\) for every nonzero \(h\in Z\).
	Since \(V\ge c>0\), the interior terms in~\eqref{eq:ball-form} are nonnegative, while the boundary terms are bounded below by \(-C_0\|h\|_{L^2(\partial B_\theta)}^2\) for some \(C_0\ge0\).
	Thus
	\[
		Q(h)\ge\frac{c}{4}\|\nabla h\|_{L^2}^2-C_0\|h\|_{L^2(\partial B_\theta)}^2.
	\]
	For every \(\varepsilon>0\), the trace inequality gives
	\[
		\|h\|_{L^2(\partial B_\theta)}^2
		\le\varepsilon\|\nabla h\|_{L^2}^2+C_\varepsilon\|h\|_{L^2}^2.
	\]
	Choosing \(\varepsilon\) sufficiently small to absorb the boundary term, we obtain constants \(\kappa>0\) and \(C<\infty\) such that
	\begin{equation}\label{eq:ball-garding}
		Q(h)\ge\kappa\|h\|_{H^1}^2-C\|h\|_{L^2}^2.
	\end{equation}
	On a bounded weakly convergent \(H^1\) sequence, the boundary traces converge strongly in \(L^2(\partial B_\theta)\) after passage to a subsequence.
	The interior terms in~\eqref{eq:ball-form} are weakly lower semicontinuous, so the same is true of \(Q\) along this subsequence.

	For a smooth test function \(\zeta\), define
	\begin{equation}\label{eq:ball-functionals}
		\begin{aligned}
			E_\zeta(h)&=\int_{B_\theta}\zeta(R_{\bar{\mathbf{g}}+h}-n(n+1))\,d\mu
			+2\int_{\partial B_\theta}\zeta(H_{\bar{\mathbf{g}}+h}-H_{\bar{\mathbf{g}}})\,d\sigma,\\
			\mathcal J(h)&=E_V(h)-\int_{\partial B_\theta} V h(\nu,\nu)(H_{\bar{\mathbf{g}}+h}-H_{\bar{\mathbf{g}}})\,d\sigma.
		\end{aligned}
	\end{equation}
	By Brendle and Marques~\cite[Theorem~10]{BrendleMarques2011}, together with \eqref{eq:ball-form} and \eqref{eq:ball-functionals}, we have
	\[
		\begin{aligned}
			|\mathcal J(h)+Q(h)|
			&\le C\int_{B_\theta}\bigl(|h||\nabla h|^2+|h|^3\bigr)\,d\mu+C\int_{\partial B_\theta}\bigl(|h|^2|\nabla h|+|h|^3\bigr)\,d\sigma.
		\end{aligned}
	\]
	The interior integral satisfies
	\[
		\begin{aligned}
			\int_{B_\theta}\bigl(|h||\nabla h|^2+|h|^3\bigr)\,d\mu
			&=\int_{B_\theta}|h|\bigl(|\nabla h|^2+|h|^2\bigr)\,d\mu\\
			&\le\|h\|_{C^0}\|h\|_{H^1}^2.
		\end{aligned}
	\]
	For the boundary integral, the \(C^1\) bound and the trace theorem give
	\[
		\begin{aligned}
			\int_{\partial B_\theta}\bigl(|h|^2|\nabla h|+|h|^3\bigr)\,d\sigma
			=&\int_{\partial B_\theta}|h|^2\bigl(|\nabla h|+|h|\bigr)\,d\sigma\\
			\le& C\|h\|_{C^1}\|h\|_{L^2(\partial B_\theta)}^2
			\le C\|h\|_{C^1}\|h\|_{H^1}^2.
		\end{aligned}
	\]
	Combining these estimates and \(\|h\|_{C^0}\le\|h\|_{C^1}\), we obtain
	\begin{equation}\label{eq:ball-first-estimate}
		|\mathcal J(h)+Q(h)|\le C\|h\|_{C^1}\|h\|_{H^1}^2.
	\end{equation}

	We next estimate \(E_\zeta(h)\) using only the \(W^{2,p}\) control of \(h\).
	Let \(\dot H_{\bar{\mathbf g}}(h)=\left.\frac{d}{dt}\right|_{t=0}H_{\bar{\mathbf g}+th}\) denote the first variation of the boundary mean curvature.
	Since \(\operatorname{div}_{\bar{\mathbf g}}h=0\) and \(h^T=0\), we have
	\begin{equation}\label{eq:boundary-mean-variation}
		2\dot H_{\bar{\mathbf g}}(h)
		=\nabla_\nu(\operatorname{tr}h)
		-\operatorname{div}_{\partial B_\theta}\bigl(h(\nu,\cdot)^T\bigr).
	\end{equation}
	The coordinate formulas for scalar curvature and mean curvature give
	\begin{equation}\label{eq:ball-remainder-bounds}
		\begin{aligned}
			\bigl|R_{\bar{\mathbf g}+h}-n(n+1)-L_{\bar{\mathbf g}}h\bigr|
			&\le C\bigl(|h||\nabla^2h|+|\nabla h|^2+|h|^2\bigr),\\
			\bigl|H_{\bar{\mathbf g}+h}-H_{\bar{\mathbf g}}-\dot H_{\bar{\mathbf g}}(h)\bigr|
			&\le C\bigl(|h||\nabla h|+|h|^2\bigr).
		\end{aligned}
	\end{equation}
	H\"older's inequality and \(\frac{p}{p-1}<2\) give
	\[
		\int_{B_\theta}|h||\nabla^2h|\,d\mu
		\le C\|h\|_{L^2}\|h\|_{W^{2,p}}.
	\]
	The integrals of the remaining interior remainder terms are bounded by \(C\|h\|_{C^1}\|h\|_{H^1}\), and the boundary integral of the absolute value of the mean curvature remainder is bounded by \(C\|h\|_{C^1}\|h\|_{L^2(\partial B_\theta)}\).
	To identify the linear terms, we combine \eqref{def:L}
	and \eqref{eq:adjoint} and integrate by parts, using \(\operatorname{div}_{\bar{\mathbf g}}h=0\):
	\begin{equation}\label{eq:ball-green-identity}
		\begin{aligned}
			&\int_{B_\theta}\bigl(\zeta L_{\bar{\mathbf g}}h
			-\langle L_{\bar{\mathbf g}}^*\zeta,h\rangle\bigr)\,d\mu\\
			=&\int_{\partial B_\theta}
			\bigl[(\operatorname{tr}h)\partial_\nu\zeta
			-\zeta\partial_\nu(\operatorname{tr}h)-h(\nu,\nabla\zeta)\bigr]\,d\sigma\\
			=&-\int_{\partial B_\theta}\zeta
			\left[\partial_\nu(\operatorname{tr}h)
			-\operatorname{div}_{\partial B_\theta}\bigl(h(\nu,\cdot)^T\bigr)\right]\,d\sigma\\
			=&-2\int_{\partial B_\theta}\zeta\dot H_{\bar{\mathbf g}}(h)\,d\sigma.
		\end{aligned}
	\end{equation}
	For the second equality, \(h^T=0\) gives \(\operatorname{tr}h=h(\nu,\nu)\) on the boundary.
	The normal derivative terms of \(\zeta\) therefore cancel, and the remaining term \(-h(\nu,\nabla_{\partial B_\theta}\zeta)\) is integrated by parts on the closed boundary.
	The last equality follows from the boundary mean curvature variation formula~\eqref{eq:boundary-mean-variation}.
	Combining \eqref{eq:ball-green-identity} with the integrated remainder bounds obtained from \eqref{eq:ball-remainder-bounds}, the Sobolev embedding \(W^{2,p}\hookrightarrow C^1\), and the trace theorem gives
	\begin{equation}\label{eq:ball-second-estimate}
		\left|E_\zeta(h)-\int_{B_\theta}\langle L_{\bar{\mathbf g}}^*\zeta,h\rangle\,d\mu\right|
		\le C_\zeta\|h\|_{W^{2,p}}\|h\|_{H^1}.
	\end{equation}

	Suppose, for contradiction, that there exists a sequence of metrics \(\mathbf g_j\) satisfying the curvature and boundary conditions, with \(\|\mathbf g_j-\bar{\mathbf g}\|_{W^{2,p}}\to0\), but none admits a boundary-fixing isometry to \(\bar{\mathbf g}\).
	After discarding finitely many terms, Lemma~\ref{lem:ball-gauge} gives diffeomorphisms \(\varphi_j\) fixing the boundary pointwise, with corresponding tensors \(h_j\) defined by~\eqref{eq:ball-perturbation}.
	These tensors are nonzero and satisfy \(\|h_j\|_{W^{2,p}}\to0\).
	Set \(\widehat h_j=\frac{h_j}{\|h_j\|_{H^1}}\).
	Passing to a subsequence gives \(\widehat h_j\rightharpoonup \widehat h\) in \(H^1\), strongly in \(L^2(B_\theta)\) and in \(L^2(\partial B_\theta)\), with \(\operatorname{div}\widehat h=0\) and \(\widehat h^T=0\).
	Apply~\eqref{eq:ball-second-estimate} with \(\zeta=V\), for which \(L_{\bar{\mathbf g}}^*V=0\).
	Using the curvature inequalities and \(V\ge c>0\) in the definition of \(E_\zeta\) in~\eqref{eq:ball-functionals}, we obtain
	\[
		0\le E_V(h_j)\le C\|h_j\|_{W^{2,p}}\|h_j\|_{H^1},
		\quad \text{and} \quad |E_\zeta(h_j)|\le\frac{\|\zeta\|_\infty}{c}E_V(h_j).
	\]
	Dividing~\eqref{eq:ball-second-estimate} by \(\|h_j\|_{H^1}\) therefore implies \(\int_{B_\theta}\langle L_{\bar{\mathbf{g}}}^*\zeta,\widehat h_j\rangle\,d\mu\to0\) for every smooth \(\zeta\).
	Passing to the weak limit gives the integral constraint in~\eqref{eq:ball-constraint-space}, hence \(\widehat h\in Z\).

	For large \(j\), \(|h_j(\nu,\nu)|<1\), so the boundary weight \(V(2-h_j(\nu,\nu))\) in \(\mathcal J(h_j)\) is positive. The curvature inequalities and~\eqref{eq:ball-functionals} then give \(\mathcal J(h_j)\ge0\).
	Estimate~\eqref{eq:ball-first-estimate} yields \(Q(\widehat h_j)\le C\|h_j\|_{C^1}\to0\).
	Weak lower semicontinuity gives \(Q(\widehat h)\le0\); the positivity on \(Z\) forces \(\widehat h=0\).
	Since \(\|\widehat h_j\|_{H^1}=1\) and \(\widehat h_j\to0\) in \(L^2\), equation~\eqref{eq:ball-garding} now gives the contradiction
	\[
		0<\kappa\le Q(\widehat h_j)+C\|\widehat h_j\|_{L^2}^2
		\le C\|h_j\|_{C^1}+C\|\widehat h_j\|_{L^2}^2\longrightarrow0.
	\]
	This proves rigidity under \(W^{2,p}\) smallness and hence under \(C^2\) smallness on the fixed compact ball.
	The smooth boundary-fixing diffeomorphism supplied by Lemma~\ref{lem:ball-gauge} therefore gives the required isometry.
\end{proof}

\section{Nonradial scalar-curvature prescription}\label{sec:nonradial}

Having proved Theorems~\ref{thm:radial} and~\ref{thm:balls}, we now give two extensions of the radial scalar-curvature prescription to nonradial scalar-curvature increments.

Recall from~\eqref{eq:adjoint} that the formal \(L^2\)-adjoint of the scalar curvature linearization is given by \(L_{\mathbf g}^*v=\Hess_{\mathbf g}v-(\Delta_{\mathbf g}v)\mathbf g-v\Ric_{\mathbf g}\). For a domain \(D\) in a smooth \((n+1)\)-dimensional  manifold, let \(\ker L_{\mathbf g}^*|_D\) denote the space of all \(v\in H^2_{\mathrm{loc}}(D)\) satisfying \(L_{\mathbf g}^*v=0\) in \(D\).
Its nonzero elements are called static potentials.
We use Corvino's localized scalar-curvature deformation theorem~\cite{Corvino2000}, with the quantitative estimate stated by Huang and Lee~\cite[Theorem~1]{HuangLee2024}.

\begin{lemma}[Corvino]\label{lem:corvino}
	Let \(D\) be a relatively compact domain with smooth boundary in a smooth \((n+1)\)-dimensional manifold, where \(n\ge2\), and let \(\mathbf g\) be a smooth metric on a neighborhood of \(\overline D\).
	Suppose \(\ker L_{\mathbf g}^*|_D=0\).
	Fix a compact set \(K\Subset D\) and \(0<\gamma<1\).
	There exist \(\varepsilon_0,C>0\) such that every \(\sigma\in C^\infty(D)\) supported in \(K\), with \(\|\sigma\|_{C^{0,\gamma}(\overline D)}<\varepsilon_0\), is realized by a symmetric two-tensor \(h\) whose extension by zero across \(\partial D\) is smooth, with \(\mathbf g+h\) positive definite and
	\[
		R_{\mathbf g+h}=R_{\mathbf g}+\sigma,
		\qquad \|h\|_{C^{2,\gamma}(\overline D)}\le C\|\sigma\|_{C^{0,\gamma}(\overline D)}.
	\]
	The constants depend on \(D,\mathbf g,K,\gamma\), and the norms are computed with respect to \(\mathbf g\).
\end{lemma}

For smooth \(\sigma\), the exponentially weighted construction gives a perturbation that extends smoothly by zero across \(\partial D\), as explained in~\cite[Sections~3.1.6 and~3.2]{Corvino2025}.
For the underlying weighted estimates, see~\cite[Section~5.1.1]{CorvinoPollack2011}.
The tensor \(h\) may be supported throughout \(\overline D\), even when the prescribed increment is supported in the smaller set \(K\).

We record the round static potentials for use in both nonradial constructions.
Since \(\Ric_{\bar{\mathbf g}}=n\bar{\mathbf g}\), taking the trace of \(L_{\bar{\mathbf g}}^*v=0\) in \eqref{eq:adjoint} gives
\[
	0=-n\Delta_{\bar{\mathbf g}}v-n(n+1)v.
\]
Hence \(\Delta_{\bar{\mathbf g}}v=-(n+1)v\), and substitution into \eqref{eq:adjoint} yields \(\Hess_{\bar{\mathbf g}}v=-v\bar{\mathbf g}\).
By interior elliptic regularity for \(\Delta_{\bar{\mathbf g}}v=-(n+1)v\), every weak solution \(v\in H^2_{\mathrm{loc}}(D)\) is smooth in \(D\).
On a connected open set \(U\subset\Sph^{n+1}\), put \(A(x)=v(x)x+\nabla v(x)\in\R^{n+2}\).
Here \(\nabla v=\operatorname{grad}_{\bar{\mathbf g}}v\), and \(\nabla^{\mathrm E}\) denotes the Euclidean connection on \(\R^{n+2}\).
By the Gauss formula for the unit sphere,
\[
	\nabla^{\mathrm E}_X\nabla v=\Hess_{\bar{\mathbf g}}v(X,\cdot)^\sharp-\langle X,\nabla v\rangle x.
\]
Since \(\nabla^{\mathrm E}_Xx=X\), \(X(v)=\langle X,\nabla v\rangle\), and \(\Hess_{\bar{\mathbf g}}v=-v\bar{\mathbf g}\), we obtain
\[
	\begin{aligned}
		\nabla^{\mathrm E}_XA
		=X(v)x+vX+\Hess_{\bar{\mathbf g}}v(X,\cdot)^\sharp-\langle X,\nabla v\rangle x
		=vX-vX=0.
	\end{aligned}
\]
Thus \(A\) is constant and every round static potential has the form
\begin{equation}\label{eq:height-functions}
	v(x)=\langle A,x\rangle\qquad\text{on }U.
\end{equation}
Conversely, each height function \(v_A(x)=\langle A,x\rangle\) satisfies \(\Hess_{\bar{\mathbf g}}v_A=-v_A\bar{\mathbf g}\), and hence \(L_{\bar{\mathbf g}}^*v_A=0\).

\begin{corollary}\label{cor:antipodal}
	Let \(\Omega\subset\Sph^{n+1}\) be a compact connected smooth domain with nonempty boundary, invariant under the antipodal map \(a(x)=-x\).
	For every \(\psi\in C_c^\infty(\operatorname{int}\Omega)\) with \(\psi\circ a=\psi\), there are smooth metrics \(\mathbf g_\varepsilon\), for all sufficiently small \(|\varepsilon|\), satisfying
	\[
		a^*\mathbf g_\varepsilon=\mathbf g_\varepsilon,
		\qquad R_{\mathbf g_\varepsilon}=n(n+1)+\varepsilon\psi,
		\qquad \mathbf g_\varepsilon\longrightarrow\bar{\mathbf g}\quad\text{in }C^2.
	\]
	They agree with \(\bar{\mathbf g}\) on a fixed neighborhood of \(\partial\Omega\), and \(\mathbf g_0=\bar{\mathbf g}\).
\end{corollary}

\begin{proof}
	Choose \(\tau>0\) sufficiently small that
	\[
		U=\{x\in\Omega:\operatorname{dist}_{\bar{\mathbf g}}(x,\partial\Omega)>\tau\}
	\]
	has smooth boundary, \(\overline U\subset\operatorname{int}\Omega\), and \(\operatorname{supp}\psi\subset U\).
	For such \(\tau\), moving the boundary inward along the round normal flow preserves connectedness, so \(U\) is connected.
	The antipodal map is a round isometry preserving \(\Omega\), hence \(a(U)=U\).
	Its free action gives a smooth relatively compact domain \(D=U/\{\pm1\}\subset\RP^{n+1}\), with quotient round metric \(\mathbf g_D\).
	Every static potential on \(D\) lifts to an even round static potential on \(U\).
	By \eqref{eq:height-functions}, the lift is a height function, which is odd; it must therefore vanish.
	Thus \(\ker L_{\mathbf g_D}^*|_D=0\).

	The function \(\psi\) descends to a compactly supported smooth function \(\psi_D\) on \(D\).
	For fixed \(0<\gamma<1\), Lemma~\ref{lem:corvino} applied to \(\sigma=\varepsilon\psi_D\) gives \(\mathbf g_D+h_\varepsilon\) with the prescribed scalar curvature and \(\|h_\varepsilon\|_{C^{2,\gamma}}\le C|\varepsilon|\|\psi_D\|_{C^{0,\gamma}}\).
	Lift this metric to \(U\) and extend it by \(\bar{\mathbf g}\) outside \(U\).
	The smooth zero extension in the lemma gives all the asserted properties; choose \(h_0=0\).
\end{proof}

For \(g\in\{2,4,6\}\), homogeneity of \(F\) gives \(F(-x)=F(x)\), so Corollary~\ref{cor:antipodal} applies to every tube \(\Sph^{n+1}_{F+c}\), where \(-1<c<1\).
For every nonzero, nonnegative, antipodally invariant \(\psi\in C_c^\infty(\operatorname{int}\Sph^{n+1}_{F+c})\) and sufficiently small \(\varepsilon>0\), the resulting metric has scalar curvature at least \(n(n+1)\), is not isometric to \(\bar{\mathbf g}\), and has unchanged boundary geometry, by the argument in the proof of Theorem~\ref{thm:radial}.
Choose a nonnegative smooth bump function \(\eta\) supported in a sufficiently small ball centered at a regular interior point, with \(\eta\) positive at its center.
The ball and its antipodal image can be chosen disjoint and small enough that their union does not cover the regular leaf through the center.
Then \(\psi=\eta+\eta\circ a\) is nonnegative, compactly supported, and invariant under the antipodal map \(a\).
Since \(\psi\) is positive at the center and vanishes at another point of the same regular leaf, it is nonradial.

For \(g=3\), we apply the local deformation theorem to nearby metrics of constant scalar curvature with no nonzero static potential.

\begin{proposition}\label{prop:cubic-prescription}
	Suppose \(g=3\), \(-1<c<1\), and \(0<\gamma<1\), and let \(K\) be a compact subset of \(\operatorname{int}\Sph^{n+1}_{F+c}\).
	For every \(\delta>0\), there exists \(\varepsilon_\delta>0\), depending also on \(F,c,K,\gamma\), such that every smooth function \(\psi\) supported in \(K\) with \(\|\psi\|_{C^{0,\gamma}}<\varepsilon_\delta\) is realized by a smooth metric \(\mathbf{g}\) satisfying
	\[
		R_{\mathbf{g}}=n(n+1)+\psi,\qquad
		\|\mathbf{g}-\bar{\mathbf{g}}\|_{C^{2,\gamma}}<\delta.
	\]
	The metrics agree with \(\bar{\mathbf{g}}\) on a common neighborhood of the boundary.
	If \(K\cap M_+=\varnothing\), they may also be chosen to equal the round metric near \(M_+\).
\end{proposition}

\begin{proof}
	By Cartan's classification, isoparametric hypersurfaces with three distinct principal curvatures occur in dimensions \(n=3,6,12,24\).
	The regular leaves of the corresponding isoparametric foliations are the principal orbits of the orthogonal actions of \(G=\mathrm{SO}(3),\mathrm{SU}(3),\mathrm{Sp}(3),F_4\), respectively, on \(\R^{n+2}\) \cite{Cartan1939,CecilRyan2025}.
	In each case, \(G\) is compact and connected, and its action is transitive on every regular leaf.
	Since \(G\) preserves the focal submanifold \(M_+\), it also preserves the radial coordinate \(r\) and the unit normal \(\partial_r\).
	The four orthogonal representations above are the isotropy representations of irreducible Riemannian symmetric spaces of rank two \cite{TakagiTakahashi1972}.
	They are therefore irreducible over \(\R\).

	Choose \(r_+<\theta\) such that the smaller open tube \(D=\{r<r_+\}\) contains \(K\).
	If \(K\cap M_+=\varnothing\), choose instead a regular annulus \(D=\{r_-<r<r_+\}\) containing \(K\), with \(0<r_-<r_+<\theta\).
	These domains are connected and \(G\)-invariant, with smooth boundary and closure in \(\operatorname{int}\Sph^{n+1}_{F+c}\).
	Choose a closed regular annulus \(\{\alpha\le r\le\beta\}\Subset D\).

	We use the constant-scalar-curvature family in~\eqref{eq:radial-family} constructed at the end of Section~\ref{sec:prescription}, with the annulus \(\{\alpha<r<\beta\}\) chosen above.
	Write \(\mathbf g_s=\mathbf g_{s,0}\) and denote the corresponding radial functions by \(N_s,u_s\).
	This smooth family is \(G\)-invariant, has scalar curvature \(n(n+1)\), agrees with \(\bar{\mathbf g}\) outside the annulus, and converges to it in every \(C^k\) norm as \(s\to0\).
	The construction also ensures that \(N_s\not\equiv1\) for every sufficiently small nonzero \(s\).

	To apply Lemma~\ref{lem:corvino}, we will rule out nonzero static potentials on \((D,\mathbf g_s)\).
	The argument below uses irreducibility to show that a nonzero static kernel would have dimension \(n+2\), which in turn would force \(\mathbf g_s\) to be Einstein.
	We therefore first show that \(\mathbf g_s\) is not Einstein for \(s\ne0\).
	Recall that \(S=\frac12\mathbf h_r^{-1}\mathbf h_r'\) and \(H=\operatorname{tr}S\) are the shape operator and mean curvature of the regular leaves with respect to the round metric \(\bar{\mathbf g}\) and the unit normal \(\partial_r\).
	Set \(\omega_s=\frac{2}{n-1}\log u_s\).
	We identify symmetric two-tensors on each regular leaf with self-adjoint endomorphisms by raising an index with \(\mathbf h_r\).
	The Gauss and Riccati equations give
	\[
		\mathbf h_r^{-1}\Ric_{\mathbf h_r}
		=(n-1)I+HS-S^2,
		\qquad S'=-I-S^2.
	\]
	For \(\widehat{\mathbf g}_{N_s}=N_s^2dr^2+\mathbf h_r\), the unit normal and shape operator are \(N_s^{-1}\partial_r\) and \(N_s^{-1}S\).
	Contracting the Gauss equation in tangential directions and using the normal Riccati equation yields
	\[
		\mathbf h_r^{-1}\Ric_{\widehat{\mathbf g}_{N_s}}|_{TM}
		=\mathbf h_r^{-1}\Ric_{\mathbf h_r}
		-\frac{S'+HS}{N_s^2}
		+\frac{N_s'}{N_s^3}S.
	\]
	Substituting the preceding identities for \(\Ric_{\mathbf h_r}\) and \(S'\) gives the tangential Ricci tensor.
	Moreover, the tangential Hessian of the radial function \(\omega_s\) is its normal derivative \(\frac{\omega_s'}{N_s}\) times the second fundamental form.
	Thus
	\begin{equation}\label{eq:cubic-ricci-hessian}
		\begin{aligned}
			\mathbf h_r^{-1}\Ric_{\widehat{\mathbf g}_{N_s}}|_{TM}
			&=(n-1+N_s^{-2})I
			+(1-N_s^{-2})(HS-S^2)
			+\frac{N_s'}{N_s^3}S,\\
			\mathbf h_r^{-1}\Hess_{\widehat{\mathbf g}_{N_s}}\omega_s|_{TM}
			&=\frac{\omega_s'}{N_s^2}S.
		\end{aligned}
	\end{equation}
	By the definition of \(\omega_s\), we have \(\mathbf g_s=e^{2\omega_s}\widehat{\mathbf g}_{N_s}\).
	In dimension \(n+1\), the conformal Ricci formula is given by
	\[
		\begin{aligned}
			\Ric_{e^{2\omega}\widehat{\mathbf g}}
			={}&\Ric_{\widehat{\mathbf g}}
			-(n-1)\bigl(\Hess_{\widehat{\mathbf g}}\omega-d\omega\otimes d\omega\bigr)-\bigl(\Delta_{\widehat{\mathbf g}}\omega
			+(n-1)|d\omega|_{\widehat{\mathbf g}}^2\bigr)\widehat{\mathbf g}.
		\end{aligned}
	\]
	We apply this formula with \(\widehat{\mathbf g}=\widehat{\mathbf g}_{N_s}\) and \(\omega=\omega_s\).
	Since \(\omega_s\) is radial, \(d\omega_s\otimes d\omega_s\) vanishes on tangential vectors.
	Restricting this formula to the regular leaf, raising an index with \(\mathbf h_r\), and substituting \eqref{eq:cubic-ricci-hessian} gives
	\begin{equation}\label{eq:cubic-ricci}
		\begin{aligned}
			\mathbf h_r^{-1}\Ric_{\mathbf g_s}|_{TM}
			={}&\left(n-1+N_s^{-2}
			-\Delta_{\widehat{\mathbf g}_{N_s}}\omega_s
			-(n-1)|d\omega_s|_{\widehat{\mathbf g}_{N_s}}^2\right)I\\
			&+\left(H(1-N_s^{-2})
			+\frac{N_s'}{N_s^3}
			-\frac{(n-1)\omega_s'}{N_s^2}\right)S
			-(1-N_s^{-2})S^2.
		\end{aligned}
	\end{equation}
	Since \(S\) has three distinct eigenvalues, the endomorphisms \(I,S,S^2\) are linearly independent.
	If \(\mathbf g_s\) were Einstein, then \(\Ric_{\mathbf g_s}=n\mathbf g_s\), and the left-hand side of \eqref{eq:cubic-ricci} would be a multiple of \(I\).
	The coefficient of \(S^2\) must therefore vanish, giving \(N_s=1\) throughout the annulus, contrary to \(N_s\not\equiv1\) for \(s\ne0\).

	We claim that \(\ker L_{\mathbf{g}_s}^*|_D=0\) for \(s\ne0\).
	Taking the trace of \(L_{\mathbf g_s}^*v=0\) in \eqref{eq:adjoint} and using \(R_{\mathbf g_s}=n(n+1)\), we obtain \(-n\Delta_{\mathbf g_s}v-n(n+1)v=0\).
	Substituting this identity into \eqref{eq:adjoint} yields
	\begin{equation}\label{eq:cubic-static}
		\Delta_{\mathbf g_s}v=-(n+1)v,\qquad
		\Hess_{\mathbf g_s}v=v(\Ric_{\mathbf g_s}-(n+1)\mathbf g_s).
	\end{equation}
	Since \(\mathbf g_s\) is smooth, interior elliptic regularity applied to the first equation in \eqref{eq:cubic-static} shows that every weak \(H^2_{\mathrm{loc}}(D)\) solution is smooth in \(D\).
	The second equation is a first-order linear system for \((v,dv)\) along curves, so a solution is determined by its value and differential at one point of the connected domain \(D\).
	In particular, the kernel has dimension at most \(n+2\).
	There is a connected \(G\)-invariant round annulus \(U\subset D\) with \(r>\beta\).
	By \eqref{eq:height-functions}, every static potential restricts on \(U\) to \(v(x)=\langle A,x\rangle\) for a unique \(A\in\R^{n+2}\).
	This defines a linear map \(v\mapsto A\) from \(\ker L_{\mathbf{g}_s}^*|_D\) into \(\R^{n+2}\).
	If \(A=0\), then \(v=0\) on \(U\), and the preceding uniqueness statement gives \(v=0\) throughout \(D\); thus the map is injective.
	For \(T\in G\), acting through the orthogonal representation above, the function \(v\circ T^{-1}\) is again a static potential and has height vector \(TA\), by the \(G\)-invariance of \(\mathbf g_s\) and \(U\).
	Thus the map is \(G\)-equivariant and its image is invariant.
	Irreducibility therefore implies that the kernel is either zero or has dimension \(n+2\).

	It remains to exclude the case \(\dim\ker L_{\mathbf g_s}^*|_D=n+2\).
	In this case, the uniqueness of static potentials from their values and differentials at a point implies that \(v\mapsto(v(p),dv(p))\) is an isomorphism onto \(\R\oplus T_p^*D\) for every \(p\in D\).
	In the following pointwise computation, the curvature tensor, gradients, and inner products are taken with respect to \(\mathbf g_s\).
	Fix \(p\in D\), and let \(X,Y\in T_pD\) be orthonormal eigenvectors of \(\Ric_{\mathbf g_s}\).
	Since \(v\mapsto(v(p),dv(p))\) is surjective, there exist static potentials vanishing at \(p\) whose gradients at \(p\) are \(X\) and \(Y\), respectively.
	For a static potential \(v\) with \(v(p)=0\), differentiating the Hessian equation in \eqref{eq:cubic-static} and commuting covariant derivatives gives, at \(p\),
	\[
		\begin{aligned}
			\langle\mathcal R(X,Y)\nabla v,W\rangle
			={}&X(v)\bigl(\Ric_{\mathbf g_s}(Y,W)-(n+1)\langle Y,W\rangle\bigr)\\
			&-Y(v)\bigl(\Ric_{\mathbf g_s}(X,W)-(n+1)\langle X,W\rangle\bigr)
		\end{aligned}
	\]
	for every \(W\in T_pD\). The terms involving derivatives of the Ricci tensor vanish because \(v(p)=0\).
	Applying this identity to the two potentials above, with \(W=Y\) and \(W=X\), respectively, yields
	\[
		\begin{aligned}
			\langle\mathcal R(X,Y)X,Y\rangle
			=\Ric_{\mathbf g_s}(Y,Y)-(n+1),\quad
			\langle\mathcal R(X,Y)Y,X\rangle
			=(n+1)-\Ric_{\mathbf g_s}(X,X).
		\end{aligned}
	\]
	By skew-symmetry of the curvature tensor in its last two arguments, we obtain \(\Ric_{\mathbf g_s}(X,X)=\Ric_{\mathbf g_s}(Y,Y)\).
	Thus all eigenvalues of \(\Ric_{\mathbf g_s}\) coincide at \(p\).
	Since \(p\) is arbitrary and \(R_{\mathbf g_s}=n(n+1)\), it follows that \(\Ric_{\mathbf g_s}=n\mathbf g_s\) on \(D\), contradicting the non-Einstein property.
	This proves the claim.

	Given \(\delta>0\), fix \(s\ne0\) sufficiently small that \(\|\mathbf{g}_s-\bar{\mathbf{g}}\|_{C^{2,\gamma}}<\frac{\delta}{2}\).
	Since \(\ker L_{\mathbf g_s}^*|_D=0\), Lemma~\ref{lem:corvino} applies to the fixed background \((D,\mathbf g_s)\) and gives, for every sufficiently small smooth \(\psi\) supported in \(K\), a smooth symmetric two-tensor \(h\) extending by zero outside \(D\), such that \(R_{\mathbf g_s+h}=n(n+1)+\psi\) and \(\|h\|_{C^{2,\gamma}}\le C_s\|\psi\|_{C^{0,\gamma}}\).
	Here norms may be taken with respect to \(\bar{\mathbf{g}}\), by equivalence of the norms for the fixed smooth metric \(\mathbf{g}_s\).
	Choose \(\varepsilon_\delta>0\) within the admissible range of Lemma~\ref{lem:corvino} and small enough that \(C_s\varepsilon_\delta<\frac{\delta}{2}\).
	For \(\|\psi\|_{C^{0,\gamma}}<\varepsilon_\delta\), the lemma gives a positive definite metric \(\mathbf g_s+h\) with \(\|h\|_{C^{2,\gamma}}<\frac{\delta}{2}\).
	Then \(\mathbf{g}=\mathbf{g}_s+h\) has the prescribed scalar curvature and lies in the required neighborhood of \(\bar{\mathbf{g}}\).
	Since \(\mathbf{g}_s=\bar{\mathbf{g}}\) and \(h=0\) outside \(D\), the resulting metric equals the round metric on a fixed neighborhood of \(\partial\Sph^{n+1}_{F+c}\).
	Choosing \(D\) to be an annulus also preserves the metric near \(M_+\).
\end{proof}

For \(0<c<1\), Proposition~\ref{prop:cubic-prescription} gives counterexamples with unchanged boundary geometry by choosing \(\psi\) to be nonnegative, nonzero, and sufficiently small.

The proposition does not, however, guarantee a differentiable family of solutions through \(\bar{\mathbf g}\) with prescribed scalar curvature \(n(n+1)+\varepsilon\psi\).
Indeed, suppose that \(\varepsilon\mapsto\mathbf g_\varepsilon\) is differentiable at \(\varepsilon=0\) in the \(C^2\) topology, with \(\mathbf g_0=\bar{\mathbf g}\) and \(R_{\mathbf g_\varepsilon}=n(n+1)+\varepsilon\psi\), and that the metrics agree with \(\bar{\mathbf g}\) on a fixed boundary neighborhood.
Then \(h=\left.\partial_\varepsilon\mathbf g_\varepsilon\right|_{\varepsilon=0}\) is compactly supported in the interior and satisfies \(L_{\bar{\mathbf g}}h=\psi\).
Since the height functions lie in \(\ker L_{\bar{\mathbf g}}^*\), integration by parts yields the necessary condition
\begin{equation}\label{eq:moment-condition}
	\int_{\Sph^{n+1}_{F+c}}\psi(x)x\,d\mu=0.
\end{equation}
A general nonradial function \(\psi\) need not satisfy \eqref{eq:moment-condition}.
The condition holds for antipodally invariant \(\psi\) when
\(g\in\{2,4,6\}\), and for radial \(\psi\) when \(g=3\):
in the latter case, the integral is a \(G\)-fixed vector and hence
vanishes by the irreducibility established above.

\medskip
\noindent\textbf{AI Assistance Disclosure.}
The authors used GPT-6 Astra to assist with language editing, literature searches, and checking the mathematical calculations and
proof details. The authors take full responsibility for the mathematical arguments,
references, and final content of the manuscript.

\end{document}